\documentclass[11pt]{article}

\usepackage{amsfonts}
\usepackage{amssymb,amsmath,amsthm}
\usepackage{latexsym}
\usepackage{fullpage}
\usepackage{hyperref}

\newtheorem{theorem}{Theorem}[section]

\newtheorem{question}{Question}[section]

\newtheorem{cor}[theorem]{Corollary}

\newtheorem{example}[theorem]{Example}

\theoremstyle{definition}

\newcounter{tenumerate}

\def\P{\mathbb{P}}

\newcommand{\one}{\1}
\newcommand{\deq}{\stackrel{\scriptscriptstyle\triangle}{=}}

\renewcommand{\epsilon}{\varepsilon}

\newcommand{\1}{\mathbf{1}}

\newcommand{\remove}[1]{}

\renewcommand{\leq}{\leqslant}
\renewcommand{\geq}{\geqslant}

\def\XXint#1#2#3{{\setbox0=\hbox{$#1{#2#3}{\int}$}
\vcenter{\hbox{$#2#3$}}\kern-.5\wd0}}

\title{Mixing under monotone censoring}
\author{Jian Ding \thanks{Partially supported by NSF grant DMS-1313596.} \\
University of Chicago \and Elchanan Mossel\thanks{
Supported by NSF grant DMS-1106999, NSF grant CCF 1320105  and DOD ONR grant N000141110140.} \\U.C. Berkeley
}

\begin{document}

\maketitle 
\begin{abstract}
We initiate the study of mixing times of Markov chain under monotone censoring. 
Suppose we have some Markov Chain $M$ on a state space $\Omega$ with stationary distribution $\pi$ and a monotone set $A \subset \Omega$. We consider the chain $M'$ which is the same as the chain $M$ started at some $x \in A$ except that moves of $M$ of the form $x \to y$ where $x \in A$ and $y \notin A$ are {\em censored} and replaced by the move $x \to x$. If $M$ is ergodic and $A$ is connected, the new chain converges to $\pi$ conditional on $A$. In this paper we are interested in the mixing time of the chain $M'$ in terms of properties of $M$ and $A$. Our results are based on new connections with the field of property testing. A number of open problems are presented. 

\end{abstract}

\section{Introduction}
\subsection{A motivating example}
Consider critical percolation on a square in the hexagonal lattice. Formally this is given by the probability space 
$\{0,1\}^{H_n}$ with the uniform distribution, where we denote by $H_n$ the sites in the hexagonal lattice in the square. 
It is trivial to sample a configuration from this model by sampling each hexagon independently. Let $A$ be the event of a left to right crossing (by $1$'s). It is well known, by duality, that $\P[A] =0.5$. Suppose we want to sample a configuration of $A$. One natural way to do so is by rejection sampling: sampling a random configuration and accepting it if and only if it is in $A$. A different natural way to sample is to start with a particular left to right crossing configuration and then repeatedly re-sample edge as long as the resulting configuration is in $A$. It is not hard to see that the second procedure will also converge to the uniform distribution on $A$. However, how long would 
it take to converge? 

\subsection{Monotone sampling in $\{0,1\}^n$}
We will study a more general question. Consider the partial order on $\{0, 1\}^n$ where $x\leq y$ if and only if $x_i \leq y_i$ for all $i\in [1, n]$.  We say a set $A\subset \{0, 1\}^n$ is monotone if $x \in A$ and $x\leq y$ imply that $y \in A$. For a monotone set $A$ and $x_0 \in A$, let $M_A^{x_0}$ denote the following Markov chain started at $x_0$. 
\begin{itemize}
\item Given the current stat $x$, pick a coordinate $i$ uniformly at random and re-randomize $x_i$ to obtain $y$. 
\item Let the next state of the chain be $y$ if $y \in A$. Otherwise let it be $x$.
\end{itemize}
It is trivial to verify that the chain converges to the uniform distribution on $A$ (it is clear that the chain is irreducible by monotonicity of $A$).  We aim to analyze the mixing time (see, e.g., \cite{AF, LPW09} for definition) for the chain $M_A^{x_0}$.

To this end, we will use a standard geometric bound on the mixing time given by the \emph{conductance} of the underlying graph for a Markov chain. Given a graph $G=G(V, E)$, the conductance $\phi(G)$ is defined to be
\begin{equation}\label{eq-def-conductance}
\phi(G) = \min_{S\subseteq V: |\mathrm{vol}(S)| \leq |E|} \Phi(S)\,,  \mbox{ where } \Phi(S) \deq \frac{|\partial_E S|}{\mathrm{vol}(S)}\,,\end{equation}
where $\mathrm{vol}(S)$  is the sum of degrees over vertices in $S$ and $\partial_E(S) = \{(x, y)\in E: x\in S, y\not\in S\}$
denotes the edge boundary set of $S$. In light of this, we will view $A$ as the underlying graph for the Markov chain $M_A$. Alternatively, $A$ can be seen as the induced subgraph of the hypercube $\{0, 1\}^n$ with a suitable number of self-loops added to each vertex so that the degree is $n$ for every $x\in A$. 

In what follows, we denote by $\P$ the uniform probability measure on $\{0, 1\}^n$. 


\begin{theorem} \label{thm:main} 
For any monotone set $A\subset \{0, 1\}^n$, we have
$$\phi(A) \geq \frac{\P[A]}{16n}\,.$$
\end{theorem}

Combined with standard results in the theory of Markov chains \cite{JS89, LS88} (see also \cite[Theorem 13.14]{LPW09}), Theroem~\ref{thm:main} yields the following corollary on the mixing time of $M_A$.

\begin{cor}\label{cor}
For any monotone set $A \subset \{0, 1\}^n$,
the mixing time for the chain $M_A$ satisfies
\[
\tau_{mix}(M_A) \leq 2 (\tfrac{16 n}{\P[A]})^2 \log(4 \cdot 2^n \P(A))\,.
\]
\end{cor} 

Note that this implies that the mixing time is polynomial in $n$ as long as $A$ is large (of measure at least inverse polynomial in $n$). In particular, the mixing time for our motivating example of sampling a critical percolation configuration with a left to right crossing has mixing time at most $O(n^3)$.  Our result is tight up to polynomial factors in $n$ as the following example shows: 

\begin{example}
Assume $n\geq 2m$ and let $A = \{x : x_1 = x_2 \ldots = x_m = 1\} \cup \{x : x_{m+1} = \ldots = x_{2m} = 1\}$. 
Clearly $A$ is monotone and $\P[A] = 2^{-m+1}$. Considering $\Phi(B)$ for $B = \{x : x_1 = x_2 \ldots = x_m = 1\}\subset A$,  we see that 
$\phi(A) \leq 2^{-m}$. Similarly starting from the point $(1,\ldots,1,0,\ldots,0)$ it is easy to see that the mixing time is lower bounded by the time to hit $(1,\ldots,1)$ with probability at least $1/4$, which is lower bounded by $2^{m-4}$.
\end{example} 

Our proof uses a new ingredient in the context of mixing of Markov chain, i.e., a result from the theory of property testing.  Property testing, explicitly defined in \cite{RS96}, plays a central role in probabilistically checkable proofs. 
However, it was extended and extensively studied on its own right for checking properties such as graph properties with 
fascinating connections to many areas of combinatorics including in particular to regularity lemmas. 
It turns out that the natural algorithm which samples a number of random neighboring pairs and rejects the monotonicity if a violating pair is seen, works well for monotonicity testing \cite{GGLRS}. The key to the success of this natural testing algorithm, which is also the key to our proof of Theorem~\ref{thm:main}, is the following structural theorem on approximately monotone set. 
\begin{theorem}\cite[Theorem 2]{GGLRS} \label{thm-GGLRS}
For any set $S \subset \{0, 1\}^n$, define
\begin{align*}
\delta(S) &= (n 2^n)^{-1} | \{(x, y) \in \{0, 1\}^n \times \{0, 1\}^n: |x-y| = 1, x\leq y, x\in S, y\not\in S\}|\,,\\
\epsilon(S)& = \min\{ \P(S \oplus A): A \mbox{ is monotone }\} \,.
\end{align*}
where $\oplus$ denotes the symmetric
difference of two sets.
Then we have $\delta(S) \geq \epsilon(S)/n$.
\end{theorem} 

\subsection{Proof of Theorem~\ref{thm:main}}
Recall that $\P$ is the uniform measure on $\{0, 1\}^n$.
In light of definition \eqref{eq-def-conductance}, it suffices to prove that 
\begin{equation}\label{eq-to-show-B}
\Phi(B) \geq \frac{\P(A)}{16n}\mbox{ for all }B\subset A\mbox{ such that } \P(B)\leq \P(A)/2\,.
\end{equation}
It is clear that \eqref{eq-to-show-B} holds if $\P(A) \P(B) < 8 \cdot 2^{-n}$ since in this case by connectivity of $A$ we have that 
$$\Phi(B) \geq \frac{1}{\mathrm{vol}(B)} = \frac{1}{\P(B) n 2^n} \geq \frac{\P(A)}{8 n}\,.$$
It remains to consider the case when $\P(A) \P(B) \geq 8\cdot 2^{-n}$.  Denote by $C = A \setminus B$ 
and by $\Omega$ the collection of
monotone sets in the  hypercube $\{0, 1\}^n$. We claim that
\begin{equation}\label{eq-not-increasing}
\mbox {{\bf either }} \P(B \oplus F )\geq  \frac{\P(A) \P(B)}{16},
\mbox{ for all } F\in \Omega\,, \mbox{{\bf or }}  \P(C \oplus F) \geq
\frac{\P(A) \P(B)}{16}\,, \mbox{ for all } F\in
\Omega\,.\end{equation}  Otherwise, there exist monotone sets $B'$
and $C'$ such that 
\begin{equation}\label{eq-prime-B-C}
\P(B \oplus B') < \frac{\P(A) \P(B)}{16}\mbox{ and }
\P(C\oplus C') < \frac{\P(A) \P(B)}{16}\,.\end{equation}
 In particular, we have $\P(B') \geq \P(B)/2$ and $\P(C') \geq \frac{7}{16}\P(A)$. An application of FKG
inequality \cite{FKG} gives that
$$\P(B'\cap C') \geq \P(|B'| )\cdot \P(C') \geq \frac{7}{32} \P(A) \P(B)\,.$$
Combined with \eqref{eq-prime-B-C}, it follows that
$$\P(B\cap C) \geq \P(B'\cap C') - \P(B\oplus B') - \P(C\oplus C') \geq  \frac{7}{32} \P(A) \P(B)- 2 \frac{1}{16} \P(A) \P(B)
>0\,,$$
contradicting with the fact that $B\cap C = \emptyset$. Thus, we completed verification of \eqref{eq-not-increasing}.

Without loss of generality we assume now that $\P( B \oplus F )\geq  \frac{\P(A) \P(B)}{16},
\mbox{ for all } F\in \Omega$ (if the same holds for $C$, we just apply the following analysis to $C$ in the same manner, with the observation that $\partial_E B = \partial_E C$).  By Theorem~\ref{thm-GGLRS}, we get that
$$|\Psi(B)| \geq \frac{2^n\P(A) \P(B)}{16} \mbox{ where } \Psi(B) \deq \{(x, y)\in \{0, 1\}^n \times \{0, 1\}^n:  |x-y| = 1, x\leq y, x\in B, y\not\in B\}  \,.$$
For $(x, y) \in \Psi(B)$, we have $x\in B$ and $x\leq y$, and thus $y\in A$ since $A$ is a monotone set.  Therefore,  we get $(x, y) \in \partial_E B$, yielding that $\Psi(B) \subseteq
\partial_E B$. This implies that $|\partial_E B| \geq \frac{2^n\P(A) \P(B)}{16}$. Combined with the fact that $\mathrm{vol}(B) = n 2^n \P(B)$, it completes the proof of \eqref{eq-to-show-B} and thus the proof of the theorem.

\subsection{Discussions and open problems}

It seems plausible that the bound on the mixing time obtained in Corollary~\ref{cor} is not sharp. A case of particular interest is when $\P(A) \geq 1/2$. Indeed, we ask the following open question.
\begin{question}
Suppose that there exists a constant $c>0$ such that a monotone subset $A\subset \{0, 1\}^n$ has measure $\P(A) \geq c$. Is it true that $\tau_{mix}(M_A) \leq C n \log n$, where $C>0$ is a constant depending only on $c$?
\end{question} 

In a different direction our results suggest testing non-product measures.  
For example, suppose we wish to reproduce Theorem~\ref{thm:main}
for the Ising model on some graph $G$, where we denote by $\mu$ the stationary measure. For  this to work we will need an analogue of the testing result. 
In this setup it is natural to define for a set $S \subset \{0,1\}^n$ (identifying $0$ with $-$ and $1$ with $+$) 
\begin{align*}
\delta(S) &=  \sum_{(x, y) \in \Psi(S)} \frac{\mu(x)}{n}, \mbox{ where } \Psi(S)  =  \{  (x, y) \in \{0, 1\}^n \times \{0, 1\}^n: |x-y| = 1, x\leq y, x\in S, y\not\in S\},\\
\epsilon(S)& = \min\{ \mu(S \oplus A): A \mbox{ is monotone}\} \,.
\end{align*}

We then ask 
\begin{question}
Consider the ferromagnetic Ising model on a graph $G = (V,E)$. Under what assumptions is it the case that
$\delta(S) \geq (\epsilon(S)/n)^a$ for all $S\subset \{0, 1\}^n$ and a fixed constant $a>0$? 
\end{question}

The following example suggests that some assumptions are needed. 
Consider Curie-Weiss model (Ising model on the complete graph) at low temperature (so the stationary measure admits double wells, see \cite{ENR2, ENR1}) with $n$ sites. For convenience, suppose that $n$ is even. Let $A = \{x: \sum_i^n x_i \leq n/2\}$. We claim that $\epsilon(A) \geq 1/6$. In order to see this, let $A_k = \{x: \sum_{i=1}^n x_i = k\}$ and $A'_k = \{x: \sum_{i=1}^n x_i = n-k\}$ for $k\leq n/2$. For $x\in A_k$ and $y\in A^c$, define
$$a(x, y) = \frac{\one_{y\in A'_k, y\geq x} \mu(y)}{|\{y: y\in A'_k, y\geq x\}|} \mbox{ and so } a(x, y) = \frac{\one_{y\in A'_k, y\geq x} \mu(x)}{|\{y: y\in A'_k, y\geq x\}|} \,.$$
Thus, $\sum_{y\in A^c} a(x, y) = \mu(x)$ for all $x\in A$. In addition, by symmetry for every $y\in A^c$ we have
$$\sum_{x\leq y, x\in A} a(x, y) = \mu(y)$$
(so $a(\cdot, \cdot)$ is a mass transportation from $A$ to $A^c$ with respect to measure $\mu$). Therefore, for any monotone set $B$ we have
$$\mu(B\cap A^c) \geq \sum_{x\in B\cap A} \sum_{y\in A^c} a(x, y) = \sum_{x\in B\cap A} \mu(x) = \mu(A\cap B)\,.$$
This implies that $\mu(B\cap A^c) \geq \mu(B)/2$. Combined with the simple fact that $\mu(A) = 1/2$, it follows that $$\mu(A\oplus B) \geq \max(\mu(A) - \mu(B), \mu(B)/2) \geq 1/6\,,$$
as desired. However, it is clear that 
$$\delta(S) \leq \mu(A_{n/2})\,,$$
which is exponentially small in $n$ at low temperature \cite{ENR2, ENR1}. 

It would be interesting to further study testing other properties for various non-product distributions.  

Finally, we note that the influence to the mixing time of censoring was studied in \cite{PW13}, where it was shown that the mixing can only be delayed for Glauber dynamics on monotone spin systems by censoring some updates (the censoring is prescribed without information on what is the proposed update). In \cite{Holroyd11}, an example was given to demonstrate that censoring can indeed speed up the mixing for proper coloring. This question was then studied in \cite{FK13} in much more general settings, which introduced a certain partial order on the class of stochastically monotone Markov kernels and proved that the monotonicity of Markov chains implies monotonicity of mixing times.  These results are different from ours in at least the following two senses: (1) They focus on Markov chains with the same stationary measure while our censoring will even change the state space of the Markov chain; (2) They aim at  qualitative results which ensure monotonicity for mixing times of Markov chains under consideration, while ours aims to give a quantitative bound on the mixing time for the censored Markov chain.

\end{document}